\documentclass[10pt]{amsart}
\usepackage{amsmath}
\usepackage{amssymb}
\usepackage{graphicx}   
\usepackage[usenames,dvipsnames]{color} 

\newcommand\marginal[1]{\marginpar{\raggedright\parindent=0pt\tiny #1}}

\newcommand{\ph}\marginal{}

\newcommand\fpf{\hfill{$\blacksquare$}}
\newcommand{\var}{{\rm var}}

\def\E{\textsf{E}}

\def\P{\textsf{P}}

\def\var{\textsf{var}}

\def\la{\lambda}
\newtheorem{theorem}{Theorem}
\newtheorem{lemma}[theorem]{Lemma}
\newtheorem{corollary}[theorem]{Corollary}

\def\eps{\varepsilon}

\def\be{\begin{equation}}
\def\ee{\end{equation}}

 \usepackage{enumerate}

   \newcommand{\lab}[1]{\label{#1}}     
\newcommand\eqn[1]{(\ref{#1})}

\def\blackslug{\hbox{\kern1pt\vrule height6pt width4pt depth1pt\kern1pt}}
\def\proof{\par\noindent{\bf Proof \enspace}\rm}
\def\qed{\penalty 500\hbox{\quad\blackslug}\ifmmode\else\par
                   \vskip4.5pt plus3pt minus2pt\fi}

\def\bX{{\bf X}}
\def\bY{{\bf Y}}
\def\bk{{\boldsymbol k}}
\def\bv{{\boldsymbol v}}
\def\bx{{\boldsymbol x}}
\def\bq{{\boldsymbol q}}
\def\bmu{{\boldsymbol \mu}}
\def\bsig{{\boldsymbol \sigma}}

\def\bz{{\boldsymbol \zeta}}
\def\bt{{\boldsymbol t}}
\def\bs{{\boldsymbol s}}
\def\uv{{\bf 1}}
\def\zv{{\bf 0}}

\def\E{\textsf{E}}
\newcommand{\ex}{{\E}}

\def\P{\textsf{P}}

\def\var{\textsf{var}}

\def\la{\lambda}

\def\eps{\varepsilon}

\author[Hitczenko]
{Pawe{\l} Hitczenko${}^\dagger$ }
\thanks{$\dagger$ This research was initiated while the first author visited Monash University. He would like to thank the School of Mathematics for their support and hospitality.}
\address{Department of Mathematics, Drexel University, Philadelphia, 
PA  19104, USA} 
\email{phitczenko@math.drexel.edu}
\author[Wormald]
{Nick Wormald*} 
\thanks{* Research supported by the Australian Laureate Fellowships grant FL120100125.}
\address{School of Mathematics, Monash University, 
VIC 3800, 
Australia} 
\email{nick.wormald@monash.edu}
\title[Multivariate asymptotic normality]{Multivariate asymptotic  normality \\ determined by high moments
}


\begin{document}
\maketitle
\begin{abstract} We extend a general result showing that the asymptotic behavior of high moments, factorial or standard, of random variables, determines
the  asymptotically normality, from the one dimensional to the multidimensional setting. This approach differs from the usual moment method which requires that  the moments of each fixed order converge. We illustrate our results by considering a joint distribution of the numbers of bins (having the same, finite, capacity) containing a prescribed number of balls in a classical allocation scheme.

\end{abstract}

\section{Introduction}
It was shown in \cite{GW} that the asymptotic behavior  of  high factorial moments can be used to establish the asymptotic normality of sequences of random variables. More precisely, it was shown that if $(X_n)$ is a sequence of non--negative random variables whose factorial moments $\E[X_n]_k$ satisfy 
\[\E[X_n]_k\sim\mu_n^k\exp\left(-\frac{k^2}{2\mu_n}\left(1-\frac{\sigma_n^2}{\mu_n}\right)\right),
\]
uniformly for $c_1\mu_n/\sigma_n\le k\le c_2\mu_n/\sigma_n$ then the sequence $(X_n-\mu_n)/\sigma_n$ converges in distribution to a standard normal random variable provided $\mu_n\to\infty$ and $\mu_n=o(\sigma_n^3)$ as $n\to\infty$. This differs from the method of employing factorial moments to  prove  convergence to a Poisson distribution, which only applies for bounded $\mu_n$. In combinatorial applications, the factorial moments are usually the easiest to compute or estimate. The above result using higher moments has been used on several occasions.  (See~\cite{CD,Jane,GWang, YLFS}.) 

A corresponding result for standard moments was established in \cite{GW} as well. Again, only the leading asymptotic term for $\E X_n^k$ is required, which, if only available for $k$ fixed, gives little information on the distribution. For instance if $X_n$ is essentially normal with $\mu_n\to\infty$ and the variance is $o(\mu_n^2)$, then typically $\E X_n^k\sim\mu_n^k$, from which one cannot deduce a nontrivial distribution result. Computing correction terms in the asymptotics can be used to remedy this situation, but this can involve tedious or even daunting computations.

Our aim here is to prove a multivariate version of the Central Limit Theorem (CLT) via the convergence of high moments (see Theorems~\ref{t:standard} and~\ref{t:factorial} below), using suitable modifications and extensions of the proof of the single-variable case in~\cite{GW}.  Our results extend to cases where the covariance matrix of the `limiting' multivariate normal distribution is asymptotically singular, up to a point, and this causes complications in their statements.

 We illustrate our result by considering a joint distribution of the number of urns, each containing a prescribed numbers of balls, say $m_i$, and each  having a finite capacity $C$  in a uniform allocation model. 
 Distributions of balls into urns are a classical topic in probability theory and have been studied extensively.  We refer to \cite{JK,KSC,KB,M} for a complete account of the theory and discussions of different models, and to e.~g.  \cite{GHP,GHW, HJ} for some of the more recent developments.  
 As one example, Timashev \cite{t} considered  urns containing the number of balls with values in a given subset $A\subset \Bbb N$. Under additional assumptions on the set $A$  he established the convergence to a normal or Poisson random variable for the number of such urns containing a given number of balls. This would essentially correspond to taking the sum of the components of our multivariate vector over all $i$ such that $m_i\in A$.

Throughout, we will use boldfaced symbols to denote the $r$--dimensional vectors, where $r$ is fixed. Operations on such vectors are meant to be coordinatewise. Thus, for example, $\bv/\bx$ denotes a vector with coordinates $v_i/x_i$, $1\le i\le r$ where $v_i$ and $x_i$ are the coordinates of $\bv$ and $\bx$, respectively. Similarly, $\log\bx$ is a vector whose entries are $\log x_i$, $1\le i\le r$. We let $[x]_k=x(x-1)\dots(x-(k-1))$ denote the $k$th falling factorial of $x$ and for vectors $\bx$ and $\bk$ we let $[\bx]_\bk$ and $\bx^\bk$ denote $\prod_{i=1}^r[x_i]_{k_i}$ and $\prod_{i=1}^rx_i^{k_i}$, respectively. We denote transpose  by $'$. Recall that each positive definite real symmetric matrix $\Sigma$ has a symmetric principal square root  with positive eigenvalues, which we denote by $\Sigma^{1/2}$.

We give our main two results in the next two sections, and discuss some applications in Section~\ref{s:app}. The applications involve distributions that are asymptotically degenerate (see Sections~\ref{s:Ex23}--\ref{s:Ex234}); some of the complexity of our main results is due to the desire to transform coordinates to obtain non-degenerate distributions. 
\medskip

\noindent
{\bf Acknowledgements}
\smallskip

\noindent
 A  result similar to  our Theorem~\ref{t:factorial} was proved independently  and simultaneously  by Berzunza Ojeda,  Holmgren and Janson \cite{bhj} where it was used to study the  asymptotics of the number of fringe trees in several tree models.  In subsequent sections we include  more technical comments on the relations between the two results. Here we would like to thank the authors of \cite{bhj} for comments and discussions concerning the topics of these papers and a suggestion to help fix a minor glitch in a proof. In particular, seeing Theorem~A.1 of~\cite{bhj} led the authors of the present paper to the simplified special case in Corollary~\ref{c:fixedG}.

\section{Standard moments}
As in \cite{GW} we begin with a statement for the standard moments. Note that in the statement below the symmetry assumption on $\Sigma_n$  is clearly not a strong restriction on applicability of the theorem, since the quadratic form in~\eqn{standard} can always be written with a symmetric matrix.

\begin{theorem}\lab{t:standard}
Let $\bX_n=(X_{n,i})$, $1\le i\le r$, $n\ge1$ be random vectors,   with $X_{n,i}\ge0$. Let $\mu_{n,i}=\E X_{n,i}$ and let $\bmu_n=(\mu_{n,i})$  be the expected value vector.  
Assume that the following hold for some   $r\times r$  real  symmetric positive definite matrices $\Sigma_n$. 
 \begin{enumerate}[(i)]
\item
For some constants $c_1<c_2$,   
\begin{equation}\lab{standard}
 \E  \bX_n ^{\bk_n}\sim \bmu_{n}^{\bk_n}\exp\left(\frac12 \left(\frac{\bk_n}{\bmu_n}\right)'\Sigma_n  \left(\frac{\bk_n}{\bmu_n}\right)\right)
\end{equation}
 uniformly   for all $\bk_n$ in $   {\rm diag}(\bmu_n)  \Sigma_n^{-1/2}B$, where $B$ is the box $ [c_1,c_2]^r$ (as a set of column vectors). 
  \item Letting $s_{n,ij}=s_{ij}=(\Sigma_n^{1/2})_{ij}$, $\tilde s_{n,ij}=(\Sigma_n^{-1/2})_{ij}$ and   
$q_{n,i}=\max_j |s_{ij}|$ 
we have
\be\lab{qmu}
q_{n,i}=o(\mu_{n,i})
\ee
and
\be\lab{max}
 \max_j |\tilde s_{ij}|q_{n,i}^2 /\mu_{n,i}=o(1)
\ee
for all $i$.
\end{enumerate}
Then  
$$
\Sigma_n^{-1/2} (\bX_n-\bmu_n) \stackrel{d}{\to}N(\zv,{\bf I}).
$$
\end{theorem}

   \proof
   
 Define 
\begin{equation}\label{t'bz}
\bt'_n=(\bk_n/\bmu_n)'\Sigma_n^{1/2},\quad    \bz_n=\Sigma_n^{-1/2}  {\rm diag}(\bmu_n)  \log(\bX_n/\bmu_n). 
\end{equation}
Then $e^{\bt'_n\bz_n}=  \left(\frac{\bX_n}{\bmu_n}\right)^{\bk_n}$ and so, using assumption (i),
\begin{align}
\E e^{\bt'_n\bz_n}&= \E\left(\frac{\bX_n}{\bmu_n}\right)^{\bk_n}\nonumber\\ 
& \sim \exp\left(\frac12 \left(\frac{\bk_n}{\bmu_n}\right)'\Sigma_n  \left(\frac{\bk_n}{\bmu_n}\right)\right) \nonumber\\
& = e^{\frac12 \bt'_n \bt_n} \lab{stdnormal}
\end{align}
for all $\bk_n$ in the stated range.

Consider the transformation from $\bk$-space into $\bt$-space defined by $T(\bk')=(\bk/\bmu_n)'\Sigma_n^{1/2}$. Each $\bt$ in the box $[c_1+\epsilon,c_2-\epsilon]^r$, where $\epsilon=(c_2-c_1)/4$ say, lies in $T(C)$ for some unit box $C$   with corners having integer coordinates. By (ii), the length of $T(C)$ is $o(1)$ in each dimension, and hence $T(C)$ is contained in $[c_1,c_2]^r$. It follows that (i) and hence~\eqn{stdnormal} applies when $\bt'_n$ is replaced by any corner $\bs'_n$  of these two boxes, and the   formulae for all corners are asymptotically equal. Expressing the arbitrary $\bt'\in T(C)$ as a convex combination $\sum_j\alpha_{jn} {\bs'_n}^{(j)}$ of a subset of these corners, we can apply (iterations of) 
 H{\"o}lder's inequality to get
$$
\E e^{\bt'_n\bz_n}=\E\exp\Big(\sum_j\alpha_{jn}  \bs'_{n j}\bz_n\Big) 
\le\prod_{j}(\ex e^ {\bs'_{nj}\bz_n})^{\alpha_{jn}}\to  e^{\frac12 \bt'_n \bt_n} 
$$
by~\eqn{stdnormal}, as the sum of $\alpha_{jn}$ over $j$ is 1.
A corresponding lower bound can be obtained by similarly expressing  the closest (or equally closest) corner, $\bs'$, of $C$ to  $\bt'$  as a convex combination of $\bt'$ with some other corners $ {\hat\bs'_{n j}} $ of $T(C')$ where $C'$ is obtained from $C$ by subtracting 1 from each coordinate. (Again, this box is contained in $B$.) In the corresponding inequality, the product of contributions from the ${\hat\bs'_{n j}}$ multiply together to give asymptotically $ e^{\frac12 \bt'_n \bt_n (1-\beta_{n})} $, where $\beta_n$ is the coefficient of $\bt'$ in this convex combination. Thus
$$
\E e^{\bs'_n\bz_n} \le  o(1)+
e^{\frac12 \bt'_n \bt_n (1-\beta_n)} (\ex e^ {\bt' _n\bz_n})^{\beta_{n}}. 
$$
Now the coefficients $\beta_{n}$ must be bounded below by a constant depending only on the dimension $r$, and hence the matching asymptotic lower bound is obtained.
It follows that the  moment generating function of $\bz_n$ converges pointwise to that of the standard multivariate normal inside  the box $[c_1+\epsilon,c_2-\epsilon]^r$.  Hence $\bz_n\stackrel{d}{\to}N(\zv,{\bf I})$.

In particular, $\Sigma_n^{-1/2}  {\rm diag}(\bmu_n)  \log(\bX_n/\bmu_n) = O_p(\uv)$ where the subscript $p$ denotes a factor bounded in probability, and the bounding is entry-wise on the vector. 
 Hence $\log(\bX_n/\bmu_n)= {\rm diag}(\uv/\bmu_n)  \Sigma_n^{1/2} O_p(\uv) =O_p(\bq_n/\bmu_n)=o_p(\uv)$, where the last equality is by  \eqref{qmu}.    
 It follows that, setting  $\bY_n:=\bX_n/\bmu_n  - \uv$ we can   expand  $\log(\bX_n/\bmu_n)=\log(\uv+\bY_n)=  \bY_n+ O_p(\bq_n^2/\bmu_n^2)$ to obtain
 $$
\bz_n = \Sigma_n^{-1/2}  {\rm diag}(\bmu_n)  ( \bY_n+O_p(\bq_n^2/\bmu_n^2))=\Sigma_n^{-1/2} (\bX_n-\bmu_n+O_p( \bq_n^2/\bmu_n )).
$$
 It follows from \eqref{max} by swapping $i$ and $j$ and using the fact that $\Sigma_n^{-1/2}  $ is symmetric, that the error term is $o_p(\uv)$.  Thus $\bz_n =  \Sigma_n^{-1/2} (\bX_n-\bmu_n )+o_p( \uv) $.
Since $\bz_n\stackrel{d}{\to}N(\zv,\uv)$, we have
  by Slutsky's Theorem that
$$
\Sigma_n^{-1/2} (\bX_n-\bmu_n) \stackrel{d}{\to}N(\zv,\bf I). 
$$
\fpf

Note that the theorem's power may lessen when there are more variables included. In particular, it is easily possible that $k =O(\mu_n/\sigma_n)$, which was required for the 1-variable case, does not always hold for the range of $k_i$ that the theorem requires. This occurs in the examples in Section~\ref{s:app}.

 Assume now that $\Sigma_n$ is of the form
\be\label{fixedG}
\Sigma_n = {\rm diag}(\bsig_n )\Gamma {\rm diag}(\bsig_n)
\ee
 for some fixed symmetric   matrix  $\Gamma=[\gamma_{i,j}]$ and vector $\bsig_n =(\sigma_{n,1},\ldots,\sigma_{n,r})$.  

Note that for any symmetric matrix $A$, the $i$th diagonal entry of  $A^2$ is the sum of squares of the entries in row $i$ of $A$. Using this with $A=\Sigma_n^{1/2}$,  we see that the numbers $q_{n,i}$ defined in (ii) of  Theorem~\ref{t:standard} satisfy
$$
q_{n,i}=\max_j|s_{i,j}|\le\left(\sum_{j=1}^r s_{i,j}^2\right)^{1/2}=\left(\gamma_{i,i}\sigma^2_{n,i}\right)^{1/2}=O(\sigma_{n,i})
$$
for $1\le i\le r$.
Thus the condition \eqref{qmu} of Theorem~\ref{t:standard} is implied by (and actually equivalent to) 
\be\lab{sigmu}
\sigma_{n,i}=o(\mu_{n,i}).
\ee

Assume further that $\Gamma$ is invertible. Then 
$$
\Sigma_n^{-1} = {\rm diag}(\bsig_{n}^{-1})\Gamma^{-1} {\rm diag}(\bsig_{n}^{-1})
$$
where $\Gamma^{-1}$ is the inverse of $\Gamma$ and $\bsig_n^{-1}=(\sigma_{n,1}^{-1},\ldots ,\sigma_{n,r}^{-1} )$.
Hence, arguing as above, the $(i,j)$ entry $\tilde s_{ij}$ of $\Sigma_n^{-1/2} $ satisfies
\begin{equation}\label{sijtilde}
|\tilde s_{ij}|\le\left(\sum_{j=1}^r\tilde s_{i,j}^2\right)^{1/2}= O(\sigma_{n,i}^{-1})
\end{equation}
and thus the condition \eqref{max} of Theorem~\ref{t:standard} is also satisfied if \eqref{sigmu} holds.
Finally, it follows from \eqref{sijtilde} that the bound on the $i$th entry $k_{n,i}$ of ${\bk}_n$ is 
$$O(\mu_{n,i} \max_j|\tilde s_{ij}|)=O( \mu_{n,i}/\sigma_{n,i}).
$$
We thus obtain the following statement
\begin{corollary}\label{c:fixedG} Let $(\bX_n)$ be as in Theorem~\ref{t:standard} and let $\Sigma_n$ be of the form \eqref{fixedG} with invertible $\Gamma$ and $\bsig_n$ satisfying $\sigma_{n,i}=o(\mu_{n,i})$ for $i=1,\dots,r$. If \eqref{standard} is satisfied for all $\bk_n$ in the range $k_{n,i}=O(\mu_{n,i}/\sigma_{n,i})$, $i=1,\dots,r$ then the conclusion of Theorem~\ref{t:standard}  holds.
\end{corollary}

 \section{Factorial moments}

\begin{theorem}\lab{t:factorial} Let  $X_{n,i}$, $1\le i\le r$, $n\ge1$ be non-negative integer variables such that $\mu_{n,i}=\E X_{n,i}\to\infty$ for all $i$. Let $\Sigma_n$ satisfy condition (ii) of Theorem~\ref{t:standard} and assume additionally that  its entries satisfy 
\begin{equation}\lab{extra}
 (\Sigma_n^{-1/2})_{ij} =o( \mu_{n,i}^{-1/3}) 
\end{equation} 
for each $i$ and $j$. Let $B$ be as in Theorem~\ref{t:standard} and 
suppose  that,  uniformly  for all $\bk_n$'s in the set 
 $D:= {\rm diag}(\bmu_n)  \Sigma_n^{-1/2}B$,
 \begin{equation}
\lab{factorialmod}
 \E [\bX]_\bk\sim \bmu_n^{ \bk_n}
 \exp\left(\frac12 \left(\frac{\bk_n}{\bmu_n}\right)'\Big(\Sigma_n- {\rm diag}(\bmu_n)\Big)  \left(\frac{\bk_n}{\bmu_n}\right)\right).  
\end{equation} 
Suppose further that 
for all   $R\subseteq [r]$,
\begin{equation}\lab{overall}
\E  \prod_{i\in R} [X_{n,i }]_{k_{n,i }}=O( 1) \prod_{i\in R} [\mu_{n,i}]_{ k_{n,i}} 
\end{equation}
when the $k_{n,i}$ are restricted so that they lie in the domain $D$  
 projected onto the coordinates in $R$.
Then, the conclusion of Theorem~\ref{t:standard} holds.
\end{theorem}

\proof
Let $\bt'_n$ and $\bz_n$ be as in \eqref{t'bz} and 
define 
$$
Q_{n,i}=\frac{[X_{n,i}]_{k_{n,i}}}{[\mu_{n,i}]_{k_{n,i}}}, \quad  Q_n= \prod_{i=1}^r Q_{n,i}, \quad U_{n,i}=\left(\frac{X_{n,i}}{\mu_{n,i}}\right)^{k_{n,i}},\quad U_n= \prod_{i=1}^r U_{n,i}.
$$
We claim that it suffices  to  show 
\begin{equation}\lab{EandQ}
\E U_n\sim \E Q_n, 
\end{equation}
for all $\bk_n\in D':=  {\rm diag}(\bmu_n)  \Sigma_n^{-1/2}[(3c_1+c_2)/4,(c_1+3c_2)/4]^r$, a slightly  truncated version of the domain  $D$ relevant to~\eqn{factorialmod}.
  For such $\bk_n$, our assumption~\eqn{extra} implies
$
k_{n,i}=o(\mu_{n,i}^{2/3})$
 and hence
 \begin{equation}\lab{convert}
 [\mu_{n,i}]_{k_{n,i}}\sim  \mu_{n,i}^{k_{n,i}}\exp(-k_{n,i}^2/2\mu_{n,i}),
 \end{equation}
  and so
the assumption~\eqn{factorialmod} gives 
\begin{equation}\lab{niceQ}
\E Q_n \sim  e^{\frac12 \bt'_n \bt_n},
\end{equation}
with $\bt'_n$ defined as in~\eqn{t'bz}.
  Thus, from~\eqn{EandQ} we will
  be done by Theorem~\ref{t:standard}.

 To prove~\eqref{EandQ}, by  the triangle inequality,
\begin{align}\lab{e1}
\Big|\prod_{i=1}^rU_{n,i}-\prod_{i=1}^rQ_{n,i}\Big|&\le\sum_{j=1}^r\Big|\prod_{i=1}^{j-1}U_{n,i}(U_{n,j}-Q_{n,j})\prod_{i=j+1}^rQ_{n,i}\Big|
\nonumber\\&=\sum_{j=1}^r\prod_{i=1}^{j-1}U_{n,i}|U_{n,j}-Q_{n,j}|\prod_{i=j+1}^rQ_{n,i}. 
\end{align}
 Since for any non--negative random variable $W$ and $1\le i\le r$, 
\[\E WU_{n,i}\le \E W+ \E WU_{n,i}I_{U_{n,i}\ge 1}\le  \E W+ \E WQ_{n,i}I_{U_{n,i}\ge 1}\le \E W+ \E WQ_{n,i},\]
it follows that the expectation of the $j$th term in~\eqn{e1} is bounded by 
\begin{equation}\lab{e2}
\sum_{S\subset [j-1]}\E \bigg[|U_{n,j}-Q_{n,j}|\prod_{i\in S\cup [j+1,r]}Q_{n,i}\bigg].
\end{equation}

We will use the following easily verified fact~\cite[(2.8)]{GW}. For $a\ge b>k$, 
\begin{equation}\lab{simple}
k\log(a/b)\leq\log {[a]_k\over
[b]_k} 
\leq 
k\log\left(\frac{a-k}{b-k}\right)
\le
{k\log(a/b)\over 1-k/b}.
\end{equation}
Let $0<x<1$ be fixed.  Note that \eqref{extra} implies that $k_{n,i}=o(\mu_{n,i}^{2/3})$ for $k_{n,i}$'s in the stated range. Hence $Q_{n,j}$ is strictly positive for $n$ sufficiently large when $x \mu_{n,j}\le X_{n,j}$.  

With  $a=X_{n,j}$, $b=\mu_{n,j}$ and $k=k_{n,j}$,~\eqn{simple} implies 
  that for $X_{n,j}\ge\mu_{n,j}$, 
$$
\log U_{n,j}\le \log Q_{n,j}\le \left(1-\frac{ k_{n,j}}{\mu_{n,j}}\right)^{-1}\log U_{n,j}
$$
which implies
$$
 \left(1-\eps_n \right) \log Q_{n,j} \le  \log U_{n,j}\le \log Q_{n,j}
$$
where $\eps_n= \max_i  k_{n,i}/ x\mu_{n,i}$. 
Similarly,
  for  $x \mu_{n,j}\le X_{n,j}<\mu_{n,j}$, 
$$
 \log Q_{n,j}\le  \log U_{n,j}\le  \left(1-\eps_n \right) \log Q_{n,j},
$$
Since $Q_{n,j}\ge 1$ iff $X_{n,j}\ge\mu_{n,j}$  the above inequalities imply  for $Q_{n,j}\ge 1$ that
\begin{equation}\lab{bigQ}
0\le Q_{n,j}-U_{n,j} \le  Q_{n,j}-(Q_{n,j})^{1-\eps_n}\le \min\{Q_{n,j}, \eps_n Q_{n,j} \log Q_{n,j}\},
\end{equation}
and for $x \mu_{n,j}\le X_{n,j}<\mu_{n,j}$ that
\begin{equation}\lab{smallQ}
0\le U_{n,j}-Q_{n,j}\le  (Q_{n,j})^{1-\eps_n}-Q_{n,j}\le  \eps_n (Q_{n,j})^{1-\eps_n}\log(1/Q_{n,j}) =O(\eps_n). 
\end{equation}

To bound \eqref{e2} consider arbitrary $S\subset[j-1]$ and let  $\tilde S_j=S\cup[j+1,r]$. We split up the expectation corresponding to each such $S$ in~\eqn{e2}  as the sum of the contributions from four regions: 
\begin{enumerate}[(i)]
\item  $X_{n,j}<x\mu_{n,j}$. Here 
$ Q_{n,j}\le U_{n,j} <  x^{k_{n,j}}$.
Thus,   the expectand (i.e.\ the argument of the expectation) is bounded by 
 \[x^{k_{n,j}}\prod_{i\in \tilde S_j }Q_{n,i}\le \eps \prod_{i\in \tilde S_j}Q_{n,i}
\]
for any $\eps$ given in advance (using $0<x<1$ and $k_{n,j}\to \infty$).
\item   $x\mu_{n,j}\le X_{n,j}<\mu_{n,j}$. Here, by~\eqn{smallQ}, the expectand is at most
\[
O(\eps_n) \prod_{i\in \tilde S_j } Q_{n,i}. 
\]

\item  $1\le Q_{n,j}\le\eps_n^{-1/2}$. Here~\eqn{bigQ} gives the bound
\[
 \eps_n^{1/2}\log(\eps_n^{-1/2}) \prod_{i\in \tilde S_j } Q_{n,i}. 
\]

\item    $Q_{n,j}>\eps_n^{-1/2}$. Here~\eqn{bigQ} bounds the expectand  by $Q_{n,j} \prod_{i\in \tilde S_j} Q_{n,i} $.  Set  integer $k'_{n,j}= (1+\delta)k_{n,j}$ for arbitrarily small $\delta=\delta(n)$. Then
   the expectand is bounded in this case by
\[
  \eps_n^{\delta/2}Q_{n,j}^{1+\delta} 
\prod_{i\in \tilde S_j} Q_{n,i}. 
\]
 \end{enumerate}
We see that each  expectation in~\eqref{e2} is bounded by  $o(1)$ times
 $$
\E \prod_{i\in  R }Q_{n,i}
$$
for some $R\subset[r]$ (which may vary from case to case). Each of these quantities  is $O(1)$ by \eqref{overall} 
 which we apply in the fourth case with   $k_{n,j}$  replaced by  $k'_{n,j}$.
 For $\delta$ sufficiently small, these parameters remain in the domain $D$ permitted by the theorem's hypotheses: since $\bk_n\in D'$ and $\Sigma_n^{-1/2}$ is invertible, the boundary of $D'$ is interior to that of $D$.
  Now~\eqn{e1} and~\eqn{e2} give that $\E|U_n-Q_n|=o(1)$.

 This implies~\eqn{EandQ} as required, since~\eqn{niceQ} implies that $\E Q_n=\Omega(1)$. \fpf
 \medskip

If $\Sigma_n$ is as in the statement of Corollary~\ref{c:fixedG},  and \eqref{factorialmod} and \eqref{overall} in Theorem~\ref{t:factorial} hold uniformly for $k_{n,i}=O(\mu_{n,i}/\sigma_{n,i})$, then to apply this theorem we need condition \eqref{extra} to be satisfied.  But, in view of~\eqn{sijtilde}, this condition is implied by (and can be shown to be equivalent to)
$$
\mu_{n,i}=o(\sigma_{n,i}^{ 3}).
$$
These assumptions on $\bsig_n$, $\bmu_n$ and $\bk_n$ are identical to those made in \cite[Theorem~A.1]{bhj} which establishes the CLT for $\Sigma_n$ of the form \eqref{fixedG}. Thus, when $\Gamma$ is invertible, Theorem~A.1 of \cite{bhj} is equivalent to our Theorem~\ref{t:factorial} in the case that $\Sigma$ is of the form  \eqref{fixedG}. Note that Theorem~A.1 of \cite{bhj} also applies to singular $\Gamma$, when the asymptotic distribution is degenerate, however then the finer structure of the distribution is not determined. Applying Theorem~\ref{t:factorial} can then give more information, by transforming the variables to a non-degenerate distribution. Our first and third examples   in the next section are of such type. 
 \section{Some  applications}\lab{s:app}
We illustrate our result with the following example. Consider placements of $n$ distinguishable balls  in $N$ distinguishable bins so
that no bin contains more than $C$ balls. Denote by $M_n(N,C)$ 
the number of such placements and let $X_m$ be the number of
bins that contain exactly $m$ balls. We are interested in the joint distribution of several of these variables, say $X_{m_1},\ldots, X_{m_r}$. To avoid trivialities, we assume $n\le CN-1$. 
  Also, for simplicity we only consider $n/N $ bounded away from 0; the less interesting case of $n/N \to 0$ can easily be dealt with similarly.  
  \begin{lemma}\lab{l:factorial}
Define
 $
g(\la)=\sum_{k=0}^C\la^k/k!
$ 
and let 
 $W=W(\la)$ be a Poisson random variable with parameter $\la$
conditioned on being at most $C$, that is
\[\P(W=j)=\frac{e^{-\la}\la^j/j!}{\sum_{k=0}^Ce^{-\la}\la^k/k!}=\frac{\la^j}{j!g(\la)}.\]
Determine $\la_0$ as the unique solution of
\begin{equation}
\label{la}\frac{\lambda_0 g'(\lambda_0)}{g(\lambda_0)}=\frac nN.
\end{equation}
Then, provided that $c_1N<n$ and $CN-n\to\infty$ for some constant $c_1>0$, and
\begin{equation}\lab{kcond}
\sum_{i=1}^r k_i = o\big(  (N/\la_0)^{2/3} \big),
\end{equation}
we have
\[\emph{\E}\prod_{i=1}^r[X_{m_i}]_{k_i}\sim 
\left(\prod_{i=1}^r\mu_i^{k_i}
\right)
\exp\Big\{-\frac{(\sum_{i=1}^r k_i(m_i-\frac{n}N))^2}{2N\emph{\var}(W(\la_0))}-\frac{(\sum_{i=1}^rk_i)^2}{2N}
\Big\}.
\]
\end{lemma}
\noindent{\bf Note:} simple calculations show that if $n/N\to C$ then $1/\la\sim  1-n/CN$.
\proof  By considering all combinations of $k_i$ bins containing $m_i$ balls each for $i=1,\ldots , r$, we have
\be\label{f_mom}\E\prod_{i=1}^r[X_{m_i}]_{k_i}
=[N]_{\sum_{i=1}^rk_i}{n\choose m_1,m_1\dots, m_r} \frac{M_{n-\sum k_im_i}(N-\sum k_i,C)}{M_n(N,C)}
\ee
where in the multinomial symbol  above $m_i$ appears $k_i$ times, $i=1,\dots,r$.

To proceed, we will need the asymptotics of $M_n(N,C)$. Let $W_1,W_2,\dots,W_N$  be the iid copies of $W$. Then, 
\begin{align*}\P\bigg(\sum_{i=1}^NW_i=n\bigg)&=\sum_{k_1+\dots+k_N=n\atop 0\le k_i\le C, 1\le i\le N}\prod_{i=1}^N\P(W_i=k_i)
=\sum_{k_1+\dots+k_N=n\atop 0\le k_i\le C, 1\le i\le N}\prod_{i=1}^N\frac{\la^{k_i}}{k_i!g(\la)}\\
&
=\frac{\la^n}{g^N(\la)n!}\sum_{k_1+\dots+k_N=n\atop 0\le k_i\le C, 1\le i\le N}
{n\choose k_1,\dots, k_N} 
=\frac{\la^n}{g^N(\la)n!}M_n(N,C).
\end{align*} 
On the other hand, if we set  $\la =\la_0$ as defined in~\eqn{la}, then, using the hypothesis that $\la=o(N)$, we have that $N\var(W)\to\infty$, so that $(W_i)$ follows the  CLT. 
  Furthermore, since the numbers $\P(W=j)$,  $0\le j\le C$, are log-concave,  by the local limit theorem  (see e.g. \cite[Lemma~2]{b} for a convenient reference) we obtain 
\[\sqrt{N\var(W)}\P\bigg(\sum_{i=1}^NW_i=n\bigg)\sim\frac1{\sqrt{2\pi}}.\]
It follows that 
\[M_n(N,C)\sim\frac{n!g^N(\lambda_0)}{\lambda_0^n\sqrt{2\pi N\var(W)}}.\]

To estimate~\eqn{f_mom} we seek the behavior of 
\[\frac{M_{n-t}(N-s,C)}{M_n(N,C)}\]
where $t=\sum k_im_i$ and $s=\sum k_i$. 
The essential part is
 \[\frac{(n-t)!}{n!}\frac{\la_0^ng^{N-s}(\la)}{\la^{n-t}g^N(\la_0)},\]
where 
$\la=\la(s,t)$ is defined by
\begin{equation}\label{r}
\frac {n-t}{N-s} = \frac{\la g'(\la)}{g(\la)} 
\end{equation}
and in particular  $\la_0=\la(0,0)$.
Differentiating \eqref{r} with respect to  $t$ and $s$  yields
\be\label{lat}
-\frac1{N-s} = \frac{ (\la g'' +g')g -(g')^2\la}{g^2}\la_t
= G(\la)\la_t,\ee 
where we have set
\be\label{G} 
G(\la):=\frac{\var(W)}\la=\la\frac{g''(\la)}{g(\la)}+\frac{g'(\la)}{g(\la)}-\la\left(\frac{g'(\la)}{g(\la)}\right)^2. 
\ee
Similarly,
\be\label{las}
\frac{n-t}{(N-s)^2} = G(\la)\la_s.\ee

 Set
\[
F(s, t) =(n-t) \ln \la -(N-s) \ln g(\la).
\]
Then 
 \[\frac{\la_0^n/g^N(\la_0)}{\la^{n-t}/g^{n-s}(\la)} =
\exp\big(-(F(s,t)-F(0,0))\big)\]
and we will expand $F(s,t)$ as a power series at $(0,0)$. We have
\begin{align*}F_t(s,t)&=-\ln\la+(n-t)\frac{\la_t}\la-(N-s)\frac{g'(\la)}{g(\la)}\la_t\\&=-\ln\la+\frac{N-s}\la\left(\frac{n-t}{N-s}-\frac{\la
   g'(\la)}{g(\la)}\right)\la_t\\&=-\ln\la,\end{align*}
   where the last equality follows from \eqref{r}.
 Similarly,
\[F_s(s,t)=(n-t)\frac{\la_s}\la+\ln
g(\la)-(N-s)\frac{g'(\la)}{g(\la)}\la_s=\ln g(\la).\]
Further, 
\[F_{tt}=-\frac{\la_t}\la=\frac1{(N-s)\la G(\la)},\]
\[ F_{ss}=\frac{g'(\la)\la_s}{g(\la)}=\frac{\la_s}\la\cdot\frac{\la
  g'(\la)}{g(\la)}=\frac{\la_s}\la\cdot\frac{n-t}{N-s}=\frac{(n-t)^2}{(N-s)^3\la G(\la)},\]
and 
\[F_{ts}=-\frac{\la_s}\la=-\frac{n-t}{(N-s)^2\la G(\la)}.
\] 
For the third order derivatives we obtain:
\[F_{ttt}=-\frac{\la_{tt}\la-\la_t^2}{\la^2},\quad
F_{tts}=-\frac{\la_{ts}\la-\la_t\la_s}{\la^2},\quad 
F_{tss}=-\frac{\la_{ss}\la-\la_s^2}{\la^2},\]
and
\[F_{sss}=-\frac{\la_{ss}\la-\la_s^2}{\la^2}\frac{n-t}{N-s}-\frac{\la_s}\la\frac{n-t}{(N-s)^2}=
\frac{n-t}{N-s}\left(\frac{\la_{ss}}\la-\frac{\la_s}\la\left(\frac{\la_s}\la-\frac1{N-s}\right)\right).\]

To get the second order partials of $\la$ 
we differentiate \eqref{lat} and \eqref{las} with respect to
$t$ and $s$ again. This gives
\[0=G'(\la)\la_t^2+G(\la)\la_{tt}
  ,\]
  i.e.
  \[\la_{tt}=-\la_t^2\frac{G'(\la)}{G(\la)}.
  \]
  We similarly get
  \[\la_{ss}=\frac{2(n-t)}{(N-s)^3G(\la)}-\la_s^2\frac{G'(\la)}{G(\la)}
  ,\]
  and
\[\la_{st}=-\frac1{(N-s)^2G(\la)}-\la_s\la_t\frac{G'(\la)}{G(\la)}.
\]
Using \eqref{lat} and \eqref{las}  we see that
\begin{align*}
\la_{tt}=&-\frac{G'(\la)}{(N-s)^2G^3(\la)},\\
\la_{ss}=&\frac{2(n-t)}{(N-s)^3G(\la)}-\frac{(n-t)^2G'(\la)}{(N-s)^4G^3(\la)},\\
\la_{st}=&-\frac1{(N-s)^2G(\la)}+\frac{(n-t)G'(\la)}{(N-s)^3G^3(\la)}.
\end{align*}
The coefficient  of $\la^{2C-1}$  in the numerator of  $G(\la)$ is zero and the same holds for the coefficient of $\la^{2C-2}$ in the numerator of $G'(\la)$. Therefore, $G(\la)=\Theta(\la^{-2})$ and $G'(\la)=O(\la^{-3})$ as long as $\la$ is bounded below by some constant $c_0$ (which we effectively show below). 
It follows that 
\begin{align*}
|F_{ttt}|&
=O\left(\frac{\la^2}{(N-s)^2}\right)
 \\
 |F_{tts}|&=
 O\left(\frac{(n-t)\la^2}{(N-s)^3}\right)\\
 |F_{tss}|&=O\left(\frac{(n-t)^2\la^2}{(N-s)^4}\right)\\
 |F_{sss}|&=O\left(\frac{(n-t)^3\la^2}{(N-s)^5}\right).
\end{align*}

Using the assumption~\eqn{kcond}  in the lemma statement, we have 
\begin{equation}\lab{splust}
s+t=O\Big(\sum k_i\Big) = o\big(  (N/\la_0)^{2/3} \big) = o(N/\la_0)=o(N)
\end{equation}
 since $c_1<\la_0=O(N)$ (the upper bound following from the note after the lemma statement). Thus, in particular, $N-s\sim N$.  Recalling the bound $G(\la)=\Theta(\la^{-2})$ from above, it now follows from~\eqn{lat} that $\la_t=O(1/NG(\la))= O(\la^2/N)$ provided that $\la\ge c_0$. Similarly, under the same assumptions, we have the same bound on $\la_s$ from~\eqn{las}. It follows from these bounds on the derivatives and~\eqn{splust} that the assumption $\la>c_0$ is justified (perhaps the easiest way to see that the restrictions on $\la$ do not affect this bound is to integrate the equation $d\la/dt = -\la^2/N$), and that
$$
|\la-\la_0|\le (s+t)O(\la^2/N) =o(\la)
$$
and hence $\la\sim \la_0$.  

 We now deduce that all third order derivatives are $O(\la_0^2/N^2)$, and so all third order terms in the expansion of $F(s,t)$ are bounded by  

\[O\left((s+t)^3\frac{\la_0^2}{N^2}\right)\cdot
o \left(\frac{N^2}{\la_0^2}\right)\cdot O\left(\frac{\la_0^2}{N^2}\right) =o(1) 
\]
and hence
\[
F(s,t)=F(0,0)+F_t(0,0)t+F_s(0,0)s+F_{s,t}(0,0)st+F_{ss}(0,0)\frac{s^2}2+F_{tt}(0,0)\frac{t^2}2+ o(1). 
\]
We thus obtain 
\begin{align*}&\frac{\la_0^n/g^N(\la_0)}{\la^{n-t}/g^{N-s}(\la)} =
\exp\big(-(F(s,t)-F(0,0))\big)\\
&\sim
\exp\bigg\{-\bigg( -t\ln\la_0+s\ln g(\la_0)-\frac n{N^2\la_0G(\la_0)}st+\frac{t^2}{2N\la_0G(\la_0)}+\frac{n^2s^2}{2N^3\la_0G(\la_0)}\bigg) \bigg\} 
\\&\quad=\frac{\la_0^t}{g^s(\la_0)}\exp\bigg\{-\frac1{2N\la_0G(\la_0)}(t-\frac{ns}N)^2\bigg\}.
\end{align*} 
When we apply this with $t=\sum_{i=1}^rm_ik_i$, $s=\sum_{i=1}^rk_i$ we
see that 
\[\frac{\la_0^t}{g^s(\la_0)}=\prod_{i=1}^r\left(\frac{\la_0^{m_i}}{g(\la_0)}\right)^{k_i}
\]
so that 
\begin{align*}&\frac{M_{n-\sum k_im_i}(N-\sum k_i,C)}{M_n(N,C)}\\&\quad
\sim
\frac{(n-\sum
  k_im_i)!}{n!}
\prod_{i=1}^r\left(\frac{\la_0^{m_i}
}{g(\la_0)}\right)^{k_i}\exp\left\{-\frac1{2N\la_0G(\la_0)}\Big(\sum k_i(m_i-\frac{n}N)\Big)^2\right\}.
\end{align*}

Combining this  with \eqref{f_mom} we get 
\begin{align*}\E\prod_{i=1}^r[X_{m_i}]_{k_i}&\sim[N]_{\sum_{i=1}^rk_i}
\prod_{i=1}^r\left(\frac{\la_0^{m_i}
}{g(\la_0)m_i!}\right)^{k_i}
\exp\Big\{-\frac{(\sum_{i=1}^r k_i(m_i-\frac{n}N))^2}{2N\la_0G(\la_0)}\Big\}\\&
\sim
\prod_{i=1}^r\left(\frac{N\la_0^{m_i}
}{g(\la_0)m_i!}\right)^{k_i}
\exp\Big\{-\frac{(\sum_{i=1}^r k_i(m_i-\frac{n}N))^2}{2N\la_0G(\la_0)}-\frac{(\sum_{i=1}^rk_i)^2}{2N}\Big\},
\end{align*}
using  $\sum_{i=1}^rk_i=o(N^{2/3})$ from~\eqn{splust}. Finally, since 
\be\label{mui}\mu_i=\E X_{m_i}=\frac{N\la_0^{m_i}
}{g(\la_0)m_i!},
\ee
and $\la G(\la) = \var(W)$, we obtain the result stated in the lemma. \fpf

We will be applying Theorem~\ref{t:factorial} to the situation where $n=CN-o(N)$. In this case $\la_0\to\infty$.  From now on, we will use $\la$ to denote $ \la_0$ as defined in Lemma~\ref{l:factorial}.  We get the relations:
\begin{align*}\frac{n}{N}=\E W&=\frac{\la g'(\la)}{g(\la)}=C\frac{1+(C-1)/\la+[C-1]_2/\la^2+\dots}{1+C/\la+[C]_2/\la^2+\dots}\\&= C\Big(1-\frac1\la+\frac{2-C}{\la^2}\Big)+O(\la^{-3}),
\end{align*}
\be\label{var}\var W=\E [W]_2+\E W-(\E W)^2 = \frac C\la\Big(1+2\frac{C-2}\la \Big)+O(\la^{-3}),\ee
since 
\[\E [W]_2=\frac{\la^2g''(\la)}{g(\la)}
=[C]_2\frac{1+(C-2)/\la+[C-2]_2/\la^2+\dots}{1+C/\la+[C]_2/\la^2+\dots}\, .
\]
 In particular, note that the number of extra  balls required to fill all bins to their maximum capacity is $CN-n =CN/\la + O(N/\la^2)$.
  
 \newcommand{\Sig}{\sigma}

We will consider several interesting instances of $\{m_1,m_2,\ldots\}\subseteq\{C-1,C-2,\ldots\}$. 
 We reserve $\Sigma$ for the appropriate matrix in \eqref{factorialmod} in Theorem~\ref{t:factorial} in the case that $m_i=C-i$ for $1\le i\le C$ and give special names to two subcases. To identify  $\Sigma$  we consider the coefficients $K_{i,i}$ and $K_{i,j}$ of $k_i^2$ and $k_ik_j$,   respectively, in the main exponential factor in the estimate for the factorial moment given by Lemma~\ref{l:factorial}.   For this example we use $\Sig_{ij}$
 to denote the $(i,j)$ element of the matrix $\Sigma$ under discussion.  We have 
\begin{equation}\lab{sii}
\frac{\Sig_{ii}}{2\mu_i^2}=
K_{i,i}+\frac1{2\mu_i}=-\frac{(m_i-\frac nN)^2}{2N\var(W(\la_0))}-\frac1{2N}+\frac1{2\mu_i}.
\end{equation}
It follows from \eqref{var} and \eqref{mui}, respectively, that
\[\frac{(m_i-\frac nN)^2}{N\var W}=O\Big(\frac\la N\Big),\quad\frac1{\mu_i}=\frac{g(\la)m_i!}{N\la^{m_i}}=\frac{m_i!\la^{C-m_i}}N\left(\frac{g(\la)}{\la^C}\right) \sim\frac{m_i!\la^{C-m_i}}{NC!}
\]
and hence, for $m_i<C-1$, 
$K_{i,i}=o(1/\mu_i)$ and $\Sig_{ii}\sim \mu_i$.
  When $m_i=C-1$ we have 
\[\frac{(C-1-\frac nN)^2}{N\var W}\sim\frac{\la(1-C/\la)^2}{NC(1+2(C-2)/\la)}\sim\frac{\la}{NC}\Big(1-\frac{2C}\la\Big)\Big(1-\frac{2(C-2)}\la\Big)
\]
and
\[
\frac1{\mu_{C-1}}=\frac{(C-1)!\la}{N}\left(\frac{g(\la)}{\la^C}\right) =\frac{\la}{NC}+\frac1{N}+\frac{C-1}{N\la}+\cdots
\]
so that for $m_i=C-1$ we have
\[
\frac{\Sig_{ii}}{2\mu_i^2}\sim\frac12\Big(-\frac\la{NC}+\frac{4(C-1)}{N{C}}+\frac\la{NC}+\frac{1}{N}-\frac{1}{N}\Big)\sim\frac {2(C-1)}{CN}.
\]
Similarly, for $0\le i\ne j< C$, 
\begin{align}
\frac{\Sig_{ij}}{\mu_i\mu_j}   =  K_{i,j}
&= -\frac{(m_i-\frac nN)(m_i-\frac nN)}{ N \var(W(\la_0))}-\frac1{2N}\label{sij} \\ 
&=\frac{-\la}{CN}
(C-m_i)(C-m_j)\Big(1+O(1/\la)\Big)\sim\frac{-\la(C-m_i)(C-m_j)}{CN}.\nonumber 
\end{align}
These estimates, together with  $\mu_i\sim NC!\la^{m_i-C}/m_i!$, 
give  the leading term asymptotics of the entries $\Sig_{ij}$ of $\Sigma$. 

In summary, if $m_i=C-i$ for $i=1,\ldots, C$, then up to the leading asymptotic term in each entry, $\Sigma$ is  
 \begin{equation}
 N\left[\begin{array}{ccccc}
\frac{4[C]_2}{\la^2}
&-\frac{2[C]_2}{\la^2}
&\dots&-\frac{C!(C-1)}{\la^{C-1}}
&-\frac{C!C}{\la^C}\\
&&&&\\
-\frac{2[C]_2}{\la^2}
&\frac{[C]_2}{\la^2}
&\dots&-\frac{(C-1)![C]_2}{\la^{C}}&-\frac{2C![C]_2}{\la^{C+1}}
\\
\vdots&\vdots&\ddots&\vdots&\vdots\\
-\frac{C!(C-1)}{\la^{C-1}}
&-\frac{(C-1)![C]_2}{\la^{C}}&\dots&\frac{C!}{\la^{C-1}}
&-\frac{(C!)^2(C-1)}{\la^{2(C-1)}}
\\
&&&&\\
-\frac{C!C}{\la^C}&-\frac{2C![C]_2}{\la^{C+1}}
&\dots&-\frac{(C!)^2(C-1)}{\la^{2(C-1)}}
&\frac{C!}{\la^C}
\end{array}\right].
\end{equation}

 To compute $\Sigma^{-1/2}$ to the required accuracy even for two bin sizes, it turns out that we will need more refined asymptotics.  As is seen in the computation of $\Sig_{ii}$ above, 
the leading term cancels on the diagonal. There are further cancellations of the next order terms when we carry out computations of the eigenvalues and eigenvectors. Below we present  a summary of some calculations we have done (with some assistance from Maple) for a few small subsets of indices. We retain   $X_i$ as the number of bins containing $i$ balls, but for other variables that depend only on the bin capacity we have found it easier to let the index $i$ refer to bins containing  $C-i$ balls  (so  $m_i=C-i$, and $\mu_i$ refers to the expected value of $X_{C-i}$).  We do not include $X_C$ in any of these cases; our main theorem still gives results for such cases, but they are less illuminating due to the larger amount of degeneracy coming from higher dependency between $X_C$ and the other variables.

\subsection{Bins containing $C-1$ and $C-2$ balls.} \lab{s:Ex23}
Here we set $r = 2$,   $m_1=C-1$ and $m_2=C-2$. Using  one more term in the expansion  \eqref{var} 
(which is required because of cancellation in  the calculation  of the leading term of the smaller eigenvalue) results  in  (noting that the $i$th row or column relates to $X_{C-i}$ in this case)
\[
\frac{\la^2}{[C]_2N}\Sigma_{1,2}= A_{1,2}:=\left[\begin{array}{cc}
4-\frac{11C+2}{\la}+O\left(\frac{1}{\la^2}\right)&-2+\frac{10C-8}\la+O\left(\frac1{\la^2}\right)\\
&\\
-2+\frac{10C-8}\la+O\left(\frac{1}{\la^2}\right)
&1-\frac{5C-4}\la+O\left(\frac{1}{\la^2}\right)
\end{array}\right].
\]
Note the near-singularity of this matrix, which   is consistent with the smaller eigenvalue tending to 0.
Thus, if $\Sigma$ is in the form of~\eqn{fixedG}, the matrix $\Gamma$   would be singular in this case, and the distribution would consequently be degenerate.

The characteristic polynomial of the matrix $A_{1,2}$, as a function of $v$, is 
\[p_{1,2}(v)=v^2-5v+\frac{(16C-2)v+9C-18}\la+O(\la^{-2}).
\]
This leads to the eigenvalues of $\Sigma_{1,2}$:
\[\nu_1=\frac{5[C]_2N}{\la^2}(1+O (\la^{-1})),\quad \nu_2=\frac{9[C]_3N}{5\la^3}(1+O (\la^{-1}))
\]
and  the corresponding normalised eigenvectors are \[
{\bf e}_1=\frac1{\sqrt5} \left[\begin{array}{c}
2+O(1/\la)\\ \\-1+O(1/\la) 
\end{array}
\right],\quad
{\bf e}_2=\frac1{\sqrt5} \left[\begin{array}{c}
1+O(1/\la)\\ \\2+O(1/\la) 
\end{array}
\right].
\]

Diagonalization then leads to  
\be\label{Sig23mh}\Sigma_{1,2}^{-1/2}
\approx
\left[\begin{array}{cc}
\frac45\nu_1^{-1/2}+\frac{1}5\nu_2^{-1/2}&-\frac{2}5\nu_1^{-1/2}+\frac25\nu_2^{-1/2}
\\ \\
-\frac{2}5\nu_1^{-1/2}+\frac25\nu_2^{-1/2}&\frac{1}5\nu_1^{-1/2}+\frac45 \nu_2^{-1/2}
\end{array}\right]
\ee 
where each of the coefficients $\frac45$ etc.\ is correct up to $O(1/\la)$ errors.
 
 The dominant entries correspond to $\nu_2^{-1/2}=\Theta(\la^{3/2}/\sqrt N)$. Thus, condition (i) in  Theorem~\ref{t:standard} 
 only requires moments for $\bk_n$ with   $k_1$ and $k_2$ both $O(\mu_i \la^{3/2}/\sqrt N)$.
Since
\[\mu_1=\E X_{C-1}\sim \frac{CN}{\la} \mbox{ and } \mu_2=\E X_{C-2}\sim \frac{[C]_2N}{\la^2} \]
this bound is $O(\sqrt{\la N})$. Hence, such $\bk_n$  satisfies the condition~\eqn{kcond} in    
  Lemma~\ref{l:factorial} provided that $ \sqrt{\la N} = o\big(  (N/\la )^{2/3} \big)$, i.e.\ 
$$
\la=o(N^{1/7}).
$$
We continue under this assumption.
 
 Regarding condition~(ii) of Theorem~\ref{t:standard}, we have  $\tilde s_{i,j}=O(\la^{3/2}/\sqrt N)$ for $i=1,2$.  Also, adjusting the last part of the argument above to estimate  $\Sigma_{1,2}^{1/2}$, we find that it has leading terms of the form \eqref{Sig23mh} with $\nu_i^{-1/2}$ replaced by $\nu_i^{1/2}$.  This time $\nu_1^{1/2}$ is dominating and thus, $q_i=O(\sqrt N/\la)$. 
 Consequently,
  condition \eqref{qmu} in Theorem~\ref{t:standard} is satisfied when $\sqrt N/\la=o(N/\la^2)$, i.e. when $\la=o(\sqrt N)$,
 and  \eqref{max} when $(\la^{3/2}/\sqrt N)N/\la^2=o(N/\la^2)$ i.e. when $\la=o(N^{1/3})$. (For both conditions, the case $i=2$ is more restrictive.) The additional assumption specific to Theorem~\ref{t:factorial}, namely condition \eqref{extra}, is satisfied whenever $\la^{3/2}/\sqrt N=o((\la^2/N)^{1/3})$ i.e. whenever $\la=o(N^{1/5})$. Thus, for  $\la= o(N^{1/7})$, Theorem~\ref{t:factorial} establishes the joint asymptotic normality for the numbers of bins containing $C-2$ and $C-1$ balls. Since $\nu_1^{-1/2}/\nu_2^{-1/2}=\frac{3\sqrt{C-2}}{ 5 \sqrt\la} (1+O(1/\la))$,    for $\Sigma^{-1/2}$ we may use the matrix
$$
A_n:=
 \frac {\la^{3/2}  }{ 3\sqrt{5[C]_3N}}   \left[\begin{array}{cc}
1+\frac{12\sqrt{C-2}}{5 \sqrt\la}& 2 - \frac{6\sqrt{C-2}}{5 \sqrt\la}
\\ \\
 2- \frac{6\sqrt{C-2}}{5 \sqrt\la} &4+\frac{3\sqrt{C-2}}{5 \sqrt\la}
\end{array}\right],
$$ 
which is correct to $O(1/\la)$ terms in the matrix entries shown. This 
provides enough accuracy to obtain a non-degenerate representation of the joint distribution: by  Theorem~\ref{t:factorial},  
$$
A_n  (\bX_n-\bmu_n) \stackrel{d}{\to}N(\zv,{\bf I}).
$$
  
 \subsection{Bins containing up to $C-2$  balls.}\label{s:C-2}

  Here, the first row and column of $\Sigma$ are irrelevant. We refer to the same matrix $\Sigma$ when indexing  the entries, but apply our theorem just to the  relevant  submatrix.  With foresight, based perhaps on the observation that for $i\ge 2$ the maximum entry in row $i$ of $\Sigma$ is $O(N/\la^i)$ on the diagonal and the off-diagonal entries are much smaller, we restrict ourselves to $k_i=O(\sqrt{\mu_i})=O(\sqrt N/\la^{i/2})$ when $i\ge 2$. (We continue to use $k_i$  for bins with $C-i$ balls, so there is no $k_1$.)
  Note that for any $j\ne i$, we have $\Sig_{ij}=O(N/\la^{i+j-1})$. Hence the contribution to the exponent in the factorial moment formula~\eqn{factorialmod} from such off-diagonal terms is
$$
O\left(\frac{N\sqrt{\mu_i}\sqrt{\mu_j}}{\la^{i+j-1} \mu_i \mu_j }\right) = O\big(\la^{1-(i+j)/2}\big)=O(1/\sqrt \la).
$$
It follows that the asymptotic relation~\eqn{factorialmod} is preserved if we replace each such entry of $\Sigma$ with 0. The diagonal terms are all $O(1)$, so we only need their leading asymptotic terms. That is, in   Theorem~\ref{t:factorial}, for the vector $(X_{C-2}, X_{C-3}, \ldots, X_0)$, we may set 
$$\Sigma = {\rm diag}(\mu_{C-2},\ldots, \mu_0),
$$
 and the moment approximation~\eqn{factorialmod} holds for the $k_i$ under consideration. Now  
 $\Sigma^{1/2}={\rm diag}(\mu_{C-2}^{1/2},\ldots, \mu_0^{1/2})$, $\Sigma^{-1/2}={\rm diag}(\mu_{C-2}^{-1/2},\ldots, \mu_0^{-1/2})$, $q_i=\sqrt{\mu_i}$, and   hypothesis (ii) in Theorem~\ref{t:standard} is satisfied. The relation~ (\ref{extra})  simply requires $\mu_{n,i}\to\infty$, which we need anyway as it is one of the main hypotheses, and finally condition~\eqn{kcond} in Lemma~\ref{l:factorial} requires just the same thing. We deduce that the  $(X_{C-2}, X_{C-3}, \ldots, X_0)$ are   asymptotically equivalent to an independent normal Gaussian with  mean vector 
 $(\mu_{C-2}, \mu_{C-3}, \ldots, \mu_0)$ and  covariance matrix ${\rm diag}(\mu_{C-2}^{1/2},\ldots, \mu_0^{1/2})$, provided that $\mu_0\to \infty$, i.e.\ $\la = o(N^{1/C})$, and $\la\to\infty$.
  
 Note that the same argument applied to any subset of these variables (but not $X_{C-1}$) yields essentially the same result, provided only that $\mu_i\to\infty$ where $i$ is the minimum index of the variables involved.
  Of course Theorem~A.1 of~\cite{bhj} also applies here, with nonsingular $\Gamma$, to give the same result.
  
We give no proof that our theorem still achieves the same result without the ``trick" of replacing the off-diagonal elements of $\Sigma$ with 0's. This is due to the awkwardness of analysing the effect of the constraints~\eqn{qmu} and~\eqn{max} in general, which depend on the size of the entries of $\Sigma^{-1/2}$. However, we carried out the computations when applying it just to the two variables $X_{C-2}$ and $X_{C-3}$, and found that these extra restrictions are satisfied provided that $\la=o(N^{1/3})$, which is required anyway by the condition that $\ex X_{C-3}\to\infty$. Presumably the general case similarly requires no extra restriction, but this does not seem worth pursuing further.
 
\subsection{Bins containing  $C-1$, $C-2$ and $C-3$ balls.}\lab{s:Ex234}
Here, with only some of the leading terms of matrix entries displayed, 
\[\Sigma_{1,2,3}\approx\frac{[C]_2N}{\la^2}\left[\begin{array}{ccc}
4-\frac{11C+2}{\la}&-2+\frac{10C-8}{\la}&-\frac{3(C-2)}{\la}\\
&&\\
-2+\frac{10C-8}{\la}&1-\frac{5C-4}{\la}
&-\frac{6[C-1]_2}{\la^2}\\&&\\
-\frac{3(C-2)}{\la}&-\frac{6[C-1]_2}{\la^2}&\frac1\la
\end{array}\right];
\]
to obtain the desired accuracy for computations we found it necessary to include further terms (not shown) correct to within $O(\la^{-4})$.
The 
 characteristic polynomial of the matrix (without leading factor) is 
 \[p_{1,2,3}(v)=v^3-5v^2+\frac{ 15Cv^2+14(C-2) v}\la+O(\la^{-2})
\]
and hence the eigenvalues of $ \Sigma_{1,2,3}$ can be expanded (to our required precision) as 
\[\nu_1\approx\frac{[C]_2N}{\la^2}\left(5-\frac{89C-28}{5\la} \right)
,\quad \nu_2\approx\frac{14[C]_3N}{5\la^3}  
,\quad \nu_3\approx\frac{8[C]_4N}{7\la^4} 
\]
with the next terms of order $1/\la$ as large as those shown.
It can then be seen that the orthogonalizing matrix entries are all $O(1)$.
\newcommand{\remove}[1]{}
\remove{
The leading terms of the orthogonalizing matrix entries are given by 
\[\left[\begin{array}{ccc}
\frac{2\sqrt5}5&\frac{3\sqrt{70}}{70}&\frac{\sqrt{14}}{14}
\\ \\
-\frac{\sqrt5}5&\frac{3\sqrt{70}}{35}&\frac{\sqrt{14}}7\\ &&\\
-\frac{6\sqrt5(C-2)}{25\la}&-\frac{\sqrt{70}}{14}&\frac{3\sqrt{14}}{14}
\end{array}\right]
\]
where the next term in $(3,1)$ entry is of order $1/\la^2$ and all others are of order $1/\la$.
}
It follows that for all $i$ and $j$ the $(i,j)$ entry $\tilde s_{i,j}$ of $\Sigma_{1,2,3}^{-1/2}$ is of order $\nu_3^{-1/2}$, i.e.\  $\tilde s_{i,j}=O(\la^2/\sqrt N)$, and that $q_i=O(\sqrt N/\la)$ for $i=1,2$ and 3.
 Since $\mu_i=\Theta(N/\la^i)$, $i=1,2,3$,
 for the relevant $\bk_n$ to satisfy the condition~\eqn{kcond} in    
  Lemma~\ref{l:factorial}, we need $\la\sqrt N = o\big((N/\la^{2/3})\big)$, i.e.\ $\la=o(N^{1/10})$. 
 
 We see that  condition  (ii) of Theorem~\ref{t:standard} holds 
  when $\la=o(N^{1/4})$ and $\la=o(N^{1/4})$, respectively, and  condition \eqn{extra} in Theorem~\ref{t:factorial} holds when $\la=o(N^{1/5})$.
Hence, we have a CLT,  
 as expressed by the conclusion of  Theorem~\ref{t:factorial}, 
 when $\la=o(N^{1/10})$. 

\subsection{Remarks on the examples}  
The covariance matrix $\Sigma_{1,2}$ is, as noted above, asymptotically singular in its leading terms, which are of comparable size. This can be explained by observing that the number of bins containing fewer than $C-2$ balls is insignificant ($\Theta(N/\la^{3})$)
 compared with $X_{C-1}=\Theta(N/\la)$ (containing $C-1$ balls), and $X_{C-2}=\Theta(N/\la^2)$. Hence $X_{C-1}$ and $X_{C-2}$ are very close to a solution to the equation obtained by counting missing balls: typically $X_{C-1}+2X_{C-2} = CN-n+O(N/\la^3)$, which explains their nearly complete dependence and  variances of comparable sizes. 

 On the other hand, the determinant of the matrix $\Sigma$  in Section~\ref{s:C-2} does not have near-cancellation in its leading term (which is of order $\la$), as the covariance is much smaller than the two variances. This can be attributed to the relevant equation being $X_{C-1}+2X_{C-2}+3X_{C-3} =   CN-n+O(N/\la^4)$. Here the variation in $X_{C-1}$ permits the covariance between $X_{C-2}$ and $X_{C-3}$ to be relatively tiny. 

Of course, $\Sigma_{1,2,3}$ contains the first example as a submatrix, and unsurprisingly given the above remarks, the main terms in its determinant cancel.   It is also not surprising that it applies for a much more restricted range of $\la$ than  the  2-variable cases. 
 
Finally, note that the main constraint in the third case, with three bin sizes, is the bound on the $k_i$  in our enumeration result  Lemma~\ref{l:factorial}. For  $k_i$ outside this bound, one would expect cubic terms to become significant in this formula, suggesting that Theorem~\ref{t:factorial} would need to be strengthened to cope with larger $\la$ in this case.

\end{document}